# Towards Prevention of Sportsmen Burnout: Formal Analysis of Sub-Optimal Tournament Scheduling

**Syed Rameez Naqvi[1], Adnan Ahmad[1], S. M. Riazul Islam[2, *], Tallha Akram[1], M. Abdullah-Al-Wadud[3] and Atif Alamri[4]**

[1]Department of Electrical and Computer Engineering, COMSATS University Islamabad, Wah Campus, Wah 47040, Pakistan
[2]Departmnet of Computer Science and Engineering, Sejong University, Seoul 05006, South Korea
[3]Department of Software Engineering, College of Computer and Information Sciences, King Saud University, Riyadh 11543, Saudi Arabia
[4]Research Chair of Pervasive and Mobile Computing, King Saud University, Riyadh 11543, Saudi Arabia
*Corresponding Author: S. M. Riazul Islam. Email: riaz@sejong.ac.kr

**Abstract:** Scheduling a sports tournament is a complex optimization problem, which requires a large number of hard constraints to satisfy. Despite the availability of several such constraints in the literature, there remains a gap since most of the new sports events pose their own unique set of requirements, and demand novel constraints. Specifically talking of the strictly time bound events, ensuring fairness between the different teams in terms of their rest days, traveling, and the number of successive games they play, becomes a difficult task to resolve, and demands attention. In this work, we present a similar situation with a recently played sports event, where a suboptimal schedule favored some of the sides more than the others. We introduce various competitive parameters to draw a fairness comparison between the sides and propose a weighting criterion to point out the sides that enjoyed this schedule more than the others. Furthermore, we use root mean squared error between an ideal schedule and the actual ones for each side to determine unfairness in the distribution of rest days across their entire schedules. The latter is crucial, since successively playing a large number of games may lead to sportsmen burnout, which must be prevented.

## 1 Introduction

Sports associations frequently organize tournaments to achieve multiple goals, such as making the league more popular and generating a specific amount of revenue. These tournaments need schedules to indicate which team will face which opponent at what time in a specified venue. Scheduling such tournaments with specific constraints is known as *Sports Scheduling* [1] in sports analytics.

Scheduling a tournament is not as simple as it looks due to many stakeholders' conflicting interests. Apart from some fundamental constraints, each stakeholder demands its own requirements. For example, television broadcaster operators require attractive matches on weekends to attract more audience [2]. In contrast, teams may request to schedule specific games on specific days due to players' non-availability [3]. Furthermore, the non-availability of venues and broadcasting equipment, due to other ongoing events, can also be a problem that should be considered during the tournament schedule [4]. Moreover, organizers always try to minimize traveling distance to decrease travel costs [5]. Constraints related to minimizing breaks and minimizing carry-over-effects are always desired to obtain a fair schedule [6,7]. Several soft and hard constraints are needed to be satisfied to mitigate such conflicts, making tournament scheduling a complex problem and is sometimes even more challenging to obtain an optimal schedule.

The problem of tournament scheduling is not new, and has been addressed a number of times by means of both research articles and scheduling tools [8,9] – this shall be covered in detail in the following section of this manuscript. In most real-life optimization problems, however, it often happens that with new events come new constraints, thereby demanding novel solutions. In outdoor sports, this is especially the case with cricket, which is currently seeing a major paradigm shift from more international to more franchise events every year. Still, however, cricket is one sport in which international events, especially those played between the rivals, always carry more weight, and therefore, attract more sponsors. This is why the International Cricket Council maintains and publishes an event calendar for the next several years in advance, and all the franchise leagues have to be scheduled in between the international events, making the former strictly time bound. As a result of this, various constraints assumed by the existing methodologies, need substantial alteration before adoption, in addition to some novel indigenous constraints demanded by each of the franchise events. For a T-20 franchise cricketing event, such as Pakistan Super League (PSL), this typically restricts (or enforces): the number of venues, weekday games, multiple games on a day, frequent traveling, more consecutive games, multiple home/away games successively, and to name a few more. While some of those may be relaxed to some extent, most of them impose hard constraints, which demand novel problem formulation and optimal schedules.

Because of the aforementioned complex constraints, it is almost impossible to find an optimal solution by manual scheduling [10]. In the past few decades, microcomputers have become powerful enough for the demanding computational tasks in the functional areas of sports scheduling. Moreover, new efficient algorithmic techniques have been developed to tackle complex computational problems. Such techniques have been employed in the past to schedule real-world sports applications. Integer programming (IP) [7,11], constraint programming (CP) [12], simulated annealing [13], branch-and-bound [14], and tabu search [15] are a few methods that have been successfully applied to solve some of the sports scheduling problems.

In this work, we focus on schedules of the previous editions of the PSL, and highlight the issues and unfairness in terms of the number of successive games and carry-over-effect, with some of the teams in particular. Our objective is to motivate for optimal scheduling of sports tournaments to avoid sportsmen burnout. The main contributions of this work are as follows:

i. We identify some of the commonly defined constraints for the tournament scheduling problem, and use them to highlight the issues and unfairness in the PSL schedule.
ii. We curate several competitive parameters for each team, including the number of home/away games, successive games, total traveling and carry-over-effect.
iii. We develop overall fatigue trends for each side using their actual schedules, and carry out a formal critical analysis of the entire schedule.
iv. We utilize cumulative distribution function (CDF) and root mean squared error as comparative measures, and propose a weighting criterion to identify the sides the PSL schedule favored the most.

Our work is important in that it establishes that the PSL schedule was tailor-made for some of the franchises, while the other sides had to suffer from increased fatigue levels due to more traveling and successive games. In one of our follow-up works (part 2 in the same series), we will propose an entire framework for optimal scheduling of the sports tournaments, specifically tailored for the T-20 leagues, which will address all the unfairness concerns with the PSL schedules. The rest of the manuscript is organized as follows: Some of the existing software for tournament scheduling and related works are reviewed in Sect. 2. Commonly defined constraints in sports scheduling, data curation and competitive parameters are presented in Sect. 3. The formal analysis and the proposed weighting criterion are presented in Sect. 4, before we conclude the manuscript in Sect. 5.

## 2 Background and Related Works

In the following, we summarize some of the related works and publicly available tools for tournament scheduling.

### 2.1 Tournament Scheduling

Tournament scheduling techniques have been used since the 1970s to solve real-world tournament scheduling problems. Campbell and Chen's work is considered an initial attempt to schedule a basketball conference of ten teams with double round-robin (DRR) constraints [8]. In the following year, Cain [16] used computer-aided heuristics to schedule Major League Baseball clubs. Many national and international sports associations have been using operation research (OR) techniques to schedule their leagues. Some of the famous leagues scheduled with scheduling methods are mentioned in the following paragraphs.

Due to Chilean geography, obtaining an optimal schedule for Chile's professional soccer leagues is a challenging and rewarding job. Chile is a very long and thin country, and for that reason, league scheduling authorities typically face plenty of constraints to make events fair and profitable. Chilean Professional Football Association (ANFP) have been using operation research techniques since 2005 to schedule Chile's professional leagues [17]. ANFP has scheduled more than 50 tournaments with OR techniques, resulting in an estimated direct economic impact of about USD 70 million, including reductions in teams traveling costs, increased ticket revenue, lower television broadcasting costs, and substantial growth in soccer television subscriptions [17].

Minor Counties Cricket Association (MCCA) league is currently scheduling their league with a schedule produced by Wright [13]. MCCA league uses a DRR format that usually completes in three years and repeats thrice over nine years for its two divisions having ten teams each. Each team needs to play three home games and three away games against its own division each year, with rotating opponents between years. Wright [13] produced a nine-year schedule with the simulated annealing technique to overcome fairness issues between teams.

### 2.2 Related Tools and Works

Various sports scheduling tools such as League Republic (LR) [9], LeagueAthletics (LA) (URL https://www.leagueathletics.com/guest/sports-scheduling-software.shtml), League Lobster (LL) (URL https://scheduler.leaguelobster.com/), All Pro League (APL) (URL https://www.allprosoftware.com/ls/), Tourney Machine (TM) (URL https://www.tourneymachine.com/Home.aspx), Print Your Brackets (PYB) (URL https://www.printyourbrackets.com/), Play Pass (PP) (URL https://playpass.com/), and Home Team Online (HTO) (URL https://www.hometeamsonline.com/), which provide plenty of options including tournament format and blackout dates to its users, are available online. A comparison between few famous tools has been conducted in Tab. 1, which justifies that tools can provide a variety of options to fulfill common constraints. However, unique constraints imposed by the leagues and excessive features available in these tools make them less acceptable, and users cannot rely on these tools to obtain an optimal schedule.

Apart from some basic constraints, each league has a set of unique constraints such as Carry-over-effect constraints and referee assignment constraints. Unfortunately, however, the available scheduling tools only provide basic constraints, which are not enough to obtain a fair schedule for a professional league. Thus, obtaining optimal schedules for professional leagues is almost impossible with these limited amounts of options available within these tools.

Moreover, tools available online – specifically APL – provide tournament scheduling option within a package, which contains various excessive features such as free league website, communication tools, and online player registration. As professional leagues usually contain these features on their own websites, they do not require these unnecessary features along with scheduling options. Thus, these tools are generally useful for small leagues and events organized by local clubs. Besides, some options that are available in one tool may not be available in another. For example, priority scheduling and algorithm selection options, available in LR, are not available in the PP scheduling tool.

While focusing on the previous works produced by researchers in the realm of cricket, this section will highlight some of the efforts that have been made to mitigate a few general issues including traveling tournament problems (TTP).

Wright introduced cost reduction objective function for England County Cricket. The latter consists of eighteen teams that compete for four different competitions, with traveling restrictions along with several other soft constraints, including team's preference constraints [18]. He used the tabu search method to reduce individual costs as well as overall cost. Wright [19] proposed a sub-cost-guided simulated annealing technique to solve a combinational optimization problem that prevailed in New Zealand Cricket. He generated several versions used in the 2003-2004 season and beyond with a tree search procedure. With a combination of IP and CP, Easton et al. [5] presented a parallel implementation of a branch-and-price algorithm to solve the TTP. While using the CP to solve feasibility problems, they used IP to solve optimality problems. Van [20] proposed a new IP based model for the Missouri Valley Conference (MVC), which organizes college basketball programs, to reduce travel costs during road trips. MVC uses travel swings that allow each team to play two consecutive matches in a single road trip. These swings, however, completely change the structure of the schedule, which affects the remaining games. Van's model assigned swings first and then placed remaining matches in empty slots. The model was flexible in assessing various scenarios by assigning swings in different ways. Duran [17] highlighted the TTP with predefined venues with integer linear programming (ILP) models. He formulated several traveling constraints to minimize the overall traveling distance covered by teams. He also used the ILP model for the Chilean first division football league to maximize the number of attractive games in final rounds. Australian Football League (AFL) consists of 18 teams that play in a single round-robin (SRR) format with five extra matches per team. As some teams have multiple home venues while others need to share the same venue, AFL encounters several fairness and cost related issues. Kyngas et al. [6] presented a three-phase process with the PEAST algorithm to mitigate break and travel distance problems in the AFL. Kostuk et al. [21] handed-down a mixed-integer programming [22] approach to develop a schedule for the Canadian Football League (CFL), providing varieties of objective functions, including maximizing the number of Friday games and minimizing the number Sunday games. They faced plenty of conflicting constraints besides structural constraints, including stadium blocks, pre-assignments, and pattern assignment constraints. However, they presented several versions to CFL for their 2010 schedule. Cavdaroglu and Atan [23] proved that rest difference problems could be solved by decomposing it into separate rounds and optimizing rounds separately. They used the IP method to minimize rest difference and developed a polynomial-time exact method for canonical schedule.

## 3 Common Constraints and Data Collection

### *3.1 Common Constraints*

There are a number of constraints which must be considered to develop a fair schedule, such as mandatory constraints of some tournament formats, round-robin for instance, and conflicting interests among stakeholders. This section will highlight some of the essential constraints that were violated in the PSL-5 schedule.

### *3.1.1 Format Constraints*

In a round-robin tournament, the most commonly used format, teams need to play against each other a fixed number of times. This format is further classified as temporarily constrained and temporarily relaxed. In a temporarily constrained round-robin tournament $r(n-1)$ slots are required; where *n* and *r* define the number of teams and number of times each team faces other teams, respectively, and the term *slot* is defined as a sufficiently long duration to accommodate a three hour game, pre- and post-match traveling, warm-up sessions, and team meetings. On a single day, there may be two such slots at most. Following decision variables are needed to formulate constraints.

$$x_{i,j,f} = \begin{cases} 1, & \text{if team } i \text{ plays home against team } j \text{ in slot } f, \\ 0, & \text{Otherwise}; \end{cases} \quad (1)$$

If set *T* has *n* number of teams and set *F* has (*n*-1) slots, the two conditions required for a single round-robin tournament are:

Each team $i \in T$ meets another $j \in T, j \neq i$, exactly once.

**Table 1:** Available Scheduling Tools

| Constraints | Tools | | | | | | | |
|---|---|---|---|---|---|---|---|---|
| | LR | APL | PYB | HTO | LL | TM | PP | LA |
| Format | ✓ | ✓ | ✓ | ✓ | ✓ | ✓ | ✓ | ✓ |
| Venue Right | ✓ | ✓ | ✗ | ✗ | ✓ | ✗ | ✗ | ✓ |
| Rest Day | ✓ | ✓ | ✗ | ✗ | ✓ | ✓ | ✗ | ✓ |
| Carry Over | ✓ | ✗ | ✗ | ✗ | ✗ | ✗ | ✗ | ✗ |
| Minimize Travel | ✗ | ✗ | ✗ | ✗ | ✗ | ✓ | ✗ | ✗ |
| Minimize Break | ✓ | ✗ | ✗ | ✗ | ✗ | ✗ | ✗ | ✓ |

$$\sum_{j \in T} \sum_{f \in F} (x_{i,j,f} + x_{j,i,f}) = 1 \qquad \forall\, i \in T,\, i \neq j \quad (2)$$

Typically used Multi Round-Robin tournaments, which are developed using mirrored Single Round-Robin, must satisfy the following constraints:

$$\sum_{i \in T} \sum_{f \in F} (x_{i,j,f} + x_{j,i,f}) = r \qquad \forall\, j \in T,\, i \neq j \quad (3)$$

$$\sum_{f \in F} (x_{i,j,f}) \geq \left\lceil \frac{r}{2} \right\rceil \qquad \forall\, i,j \in T,\, i \neq j \quad (4)$$

$$\sum_{f \in F} (x_{i,j,f} + x_{j,i,f}) = 1 \qquad \forall\, i,j \in T,\, i \neq j \quad (5)$$

Constraint (4), (5) and (6) restrict that a team *i* faces another team *j*, *r* times, number of games played at each teams venues must not differ by more than 1, and each team plays exactly once in each spot, respectively.

*3.1.2 Consecutive Break Constraints:*

Though undesirable, a specific number of breaks is required in every sports league. Two variables, *Y* and *Z*, setting the lower and upper bounds respectively can be defined to control these breaks such that:

$$Y \leq \sum_{y=0}^{Z} \sum_{j=1}^{n} (x_{i,j,f+y}) \leq Z \quad (6)$$

$$Y \leq \sum_{y=0}^{Z} \sum_{j=1}^{n} (x_{j,i,f+y}) \leq Z \quad (7)$$

Constraints (7) and (8) show that a team must play at least *Y* and at most *Z* consecutive home and away games, respectively.

*3.1.3 Break-Slot Constraints:*

Leagues often avoid breaks in second and last slot which can be achieved through following constraints:

$$\sum_{j=1}^{n} (x_{i,j,1} + x_{i,j+\hat{j},2}) = 1 \qquad \forall\, i, \hat{j} \in T, i \neq j \neq \hat{j} \quad (8)$$

$$\sum_{j=1}^{n} (x_{i,j,\hat{F}-1} + x_{i,j+\hat{j},\hat{F}}) = 1 \qquad \forall\, i, \hat{j} \in T, i \neq j \neq \hat{j} \quad (9)$$

$\hat{F}$ are the total number of slots.

*3.1.4 Traditionally Strong Team Constraint*

If consecutive games are scheduled among traditionally strong and weak teams, the latter can face issues such as fatigue or loss of interest. It can be avoided by following constraint which states that no team can play two consecutive games against strong teams:

Let $S$ be a set of strong teams.

$$\sum_{j \in S}(x_{i,j,f} + x_{j,i,f} + x_{i,j,f+1} + x_{j,i,f+1}) \leq 1 \qquad \forall\, i \in T, f \in F \tag{10}$$

*3.1.5 Complementary Constraints*

Teams that share the same venue cannot play a home match in same slot which is guaranteed by the following constraint, imposing that if two teams $i$ and $j$ share the same venue and play the games $d$ and $e$, respectively:

$$x_{i,d,f} + x_{j,e,f} \leq 1 \qquad \forall\, i,j \in T, d,e \in V, f \in F \tag{11}$$

where $V$ is a set of games.

*3.1.6 Game Constraints*

Game constraints are either used to prevent or fix certain games in certain slots based on team or broadcaster's request, player availability or popular matches, as expressed below:

Let $A_{i,j}$ be a set of player availability slots, the slots in which some players who are bound to play international sports for their country are available to play in the league. $g_{i,j}$ is a game popularity parameter with possible value 1 or 0 to indicate whether the game is attractive (1) or ordinary (0).

which is either 1 or 0 if the game is attractive or not, respectively. $g_{max}$ is introduced as upper bound limit of a match's popularity:

$$\sum_{i \in T}\sum_{j \in T} g_{i,j}\, x_{i,j,f} \leq g_{max} \qquad \forall\, f \in F \tag{12}$$

$$\sum_{p \in A_{i,j}} x_{i,j,p} = 1 \qquad \forall\, i,j \in T, i \neq j \tag{13}$$

Constraint (13) limits number of attractive games in particular slots, (14) respects player availability.

*3.1.7 Place Constraints:*

Due to unavailability of home venues, leagues impose place constraint ensuring that teams play in home venues in specific slots, defined as:

Let $B_i$ is a set of venue availability slots.

$$\sum_{p \in B_i} x_{i,j,p} = 1 \qquad \forall\, i,j \in T, i \neq j \tag{14}$$

### *3.2. Data Curation From the PSL Schedules*

Although, our model will be applicable to majority of franchise T-20 leagues, our focus in this study is PSL because of the unfairness in scheduling. The fifth edition, known as PSL-5, was to follow DRR method, where six teams participated, each representing their city, playing on four venues. Note that two teams were No-home-venue (NHV) teams, whose games were not scheduled fairly across the venues. Rather, they were assigned a different venue, with minimum aerial distance, as secondary home. For example, the teams Peshawar Zalmi (PZ) shared Rawalpindi Stadium (R) with Islamabad United (IU), Quetta Gladiators shared National Stadium (K) with Karachi Kings (KK) and Multan Sultans, despite having their home venue at Multan Stadium, shared Qaddafi Stadium (L) with Lahore Qalandars (LQ).

Typically, each team should play half of its games at home venue, called home games while rest at the opponents venue called away games. NHV teams must be treated in the same manner; however when an NHV team, $T$, plays its secondary-home game against a team $U$, such that their venue is the primary home to $U$, then the game must be termed as *away* for $T$. In such situations, it is only fair for the game to be played at a neutral venue.

*3.2.1 Number of Games in Each Venue*

Out of the total thirty round games ($2 \times C_2^6$), eleven were played in Lahore, eight were played each in the Karachi and Rawalpindi, while only three were played in the Multan. The unfairness may be due to the fact that Lahore is situated in the middle, making it an easy destination for all teams. We believe that this difference must be reduced, if not waived off. Ideally, every team deserves to play at their home grounds.

**Table 2:** Unfairness in the PSL-5 schedule w.r.t traveling

| Teams | Venue in each game | | | | | | | | | | Total Travel | Home Games |
|---|---|---|---|---|---|---|---|---|---|---|---|---|
| | 1 | 2 | 3 | 4 | 5 | 6 | 7 | 8 | 9 | 10 | | |
| IU | K | L | L | R | R | R | L | R | R | K | 5 | 5 |
| KK | K | K | M | R | R | L | L | K | K | K | 4 | 5 |
| LQ | L | L | R | L | L | L | L | L | K | L | 4 | 8 |
| MS | L | L | M | M | M | L | R | L | K | L | 6 | 3 |
| PZ | K | K | M | R | R | R | R | L | K | | 4 | 0 |
| QG | K | K | K | R | M | L | R | L | L | K | 6 | 0 |
| Legend | K: National Stadium Karachi | | | | | | | | | | | |
| | L: Gaddafi Stadium Lahore | | | | | | | | | | | |
| | R: Rawalpindi Cricket Stadium | | | | | | | | | | | |
| | M: Multan Cricket Stadium | | | | | | | | | | | |

*3.2.2 Total Travelling*

Tab. 2 summarizes the venues for each team's ten round games in order. It is evident that IU had to travel inter-cities 5 times to complete its schedule. Similarly, KK, LQ and PZ had to travel 4 times, while MS and QG, unfortunately, had to travel 6 times.

*3.2.3 Number of Home Games*

The last column of Table 2 indicates that IU, KK, LQ and MS played 5, 5, 8 and 3 three games at home venues respectively. Since, PZ and QG were NHV teams, their last columns have zeros. However, PZ and QG played 4 and 3 secondary home games respectively, which are not listed in the table.

*3.2.4 Number of Away Games*

Typically in a DRR method, if a team *T* has played an away game at *U*'s home venue, then it must play its home game its venue instead of a neutral one, which is fair. PSL-5's schedule was extremely unfair in this context, unfortunately. Such details are presented in Tab. 3 in which each row lists home games for each team. Evident from the table is that IU played its home games against KK, MS, PZ and QG in its home venue. However, its home game against LQ was unfairly played at LQ's home venue, therefore called an away game. Similarly, PZ played its four home games in its secondary home venue, including a game against IU (the primary home team for Rawalpindi Cricket Stadium), thus listed as away game. Finally, PZ played against MS in a neutral venue, despite the fact that the latter played PZ its home game in its home venue Multan Cricket Stadium. Several such unfairness examples are presented in this table; observe that LQ was the only team to play all of its home games in its home venue, while QG was the most unfortunate, playing three home games in opponents' venues, called away games.

*3.2.5 Carry-Over Effect*

Each team demands equal rest days as successive games and traveling leads to severe carryover effect on players. In time-bound events, it is almost impossible to guarantee. Thus, some teams play more games on consecutive days compared to others. Tab. 4 shows, for each team, the average number of games per day, total number of rest days between the first and the last game, and the number of games 'X' played in 'Y' consecutive days, represented by 'XinY streak'.

**Table 3:** Home-and-away map between opponents

|    | IU | KK | LQ | MS | PZ | QG |
|----|----|----|----|----|----|----|
| IU | -  | H  | A  | H  | H  | H  |
| KK | H  | -  | H  | A  | H  | H  |
| LQ | H  | H  | -  | H  | H  | H  |
| MS | S  | H  | A  | -  | H  | H  |
| PZ | A  | S  | S  | N  | -  | S  |
| QG | S  | A  | A  | A  | S  | -  |

Legend  H: Home game played at home
A: Home game played at opponent's home
S: Home game played at secondary home
N: Home game played ta a neutral venue

**Table 4:** Unfairness in the PSL-5 schedule w.r.t carry-over effect

| Team | Avg. games per day | Number of no-games days | 2in2 streak | 3in4 streak | 4in6 streak | 6in11 streak |
|------|---------|----|---|---|---|---|
| IU   | 0.416   | 14 | 3 | 2 | 0 | 2 |
| KK   | 0.416   | 14 | 2 | 2 | 0 | 0 |
| LQ   | 0.416   | 14 | 2 | 1 | 2 | 1 |
| MS   | 0.416   | 14 | 2 | 1 | 0 | 0 |
| PZ   | 0.455   | 12 | 2 | 2 | 1 | 2 |
| QG   | 0.40    | 15 | 1 | 1 | 0 | 0 |

**4 Proposed Formal Modeling and Analysis**

*4.1 Competitive Parameters*

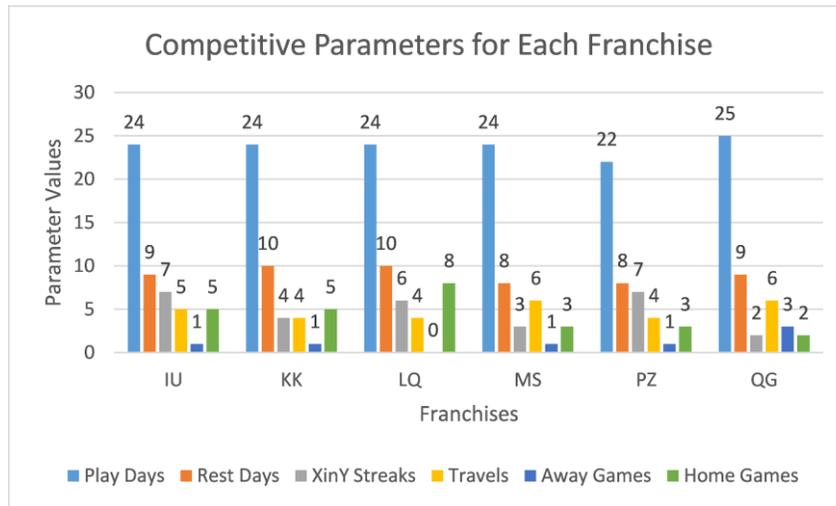

**Figure 1:** Summary of competitive parameters for each franchise

Evidently, QG is the most favored team, playing minimum number of matches on average per day, had maximum rest days, and only once playing two games in two consecutive days and three games in four consecutive days. In contrast, PZ was the most unfortunate one, playing maximum games on average per day, had minimum rest days, and played consecutive games on several occasions. The last one included PZ's two streaks of games on two consecutive days and three games in four consecutive days, one streak of four games in six consecutive days, and two streaks of six games in eleven consecutive days. This made PZ's schedule compact, resulting in less travel and more games.

Fig. 1 summarizes all the above discussed details as a histogram. The parameters plotted in the figure correspond to the following in order: 1) the total number of days that each franchise took to play all of its round games, 2) the total number of rest days (without intercity travel days), 3) the total number of XinY streaks, 4) the total number of intercity travels, 5) the total number of away games that a franchise must have played on its home/secondary venue, and 6) the total number of round games that a franchise played on its home ground.

*4.2 Formulation of Fatigue Trends*

The above histogram is self-explanatory in that it highlights how the schedule has favored some of the franchises using a number of competitive parameters. What remains unaddressed through the histogram, however, is the fact that playing successively – in long streaks – and plenty of traveling, severely affect the players' fitness, both mentally and physically. Especially if the rest days are not wisely distributed across the entire schedule for a given team, its fatigue may rise exponentially, already making it an unfair competition even before the game has commenced. For this purpose, based on the recent literature [24, 25], we have developed 30-hour general fatigue trends for each game a team plays and each time it travels using exponential utility ($U$); this is shown in Fig. 2, and given by eqn. (15).

$$\vec{U}_g/\vec{U}_t = \begin{cases} \exp(5 \times (x-1)), & \text{for the 7 hours of play/travel} \\ -(1 - \exp(-5 \times y)) + 1, & \text{for the next 23 hours of rest} \end{cases} \quad (15)$$

Here $x$ and $y$ are linearly spaced vectors between 0 and maximum fatigue values (1 for play and 0.5 for travel), and the spacing is given by $\frac{max-min}{n-1}$, where $n$ is the total number of points, e.g. 7 for play/travel (length of $x$) and 23 for rest (length of $y$).

Although, these trends are merely an estimate, in our case, they will be applied on each franchise in an identical manner, without compromising on generality. This will assist in observing the net fatigue on each franchise's players at the end of the round matches, and determine in what ways the schedule has been unfair, and where it needs improvement. Fig. 3 is the result of applying the fatigue trends on each franchise's schedule, while their equivalent histograms and cumulative distributions are given in Fig. 4 and Fig. 5 respectively. The histogram for each represents density of its fatigue regions. Clearly, IU is the most unfortunate in that it possessed shallowest (near) zero-fatigue levels – indicating the fact that they did not get enough days to rest. Similarly, the CDF in Fig. 5 validates the same: IU was the slowest to approach the maximum cumulative distribution value.

*4.3 Fatigue Analysis*

While it may be conveniently observed in Fig. 3 that some franchises played more games in quick succession, there must be some metrics for drawing a fair quantitative comparison for teams' fatigues as well. In what follows, we enumerate three of those metrics, which - to the best of our knowledge - should enable the readers to clearly see the franchise this scheduled favored the most, as far as teams' overall fatigue is concerned.

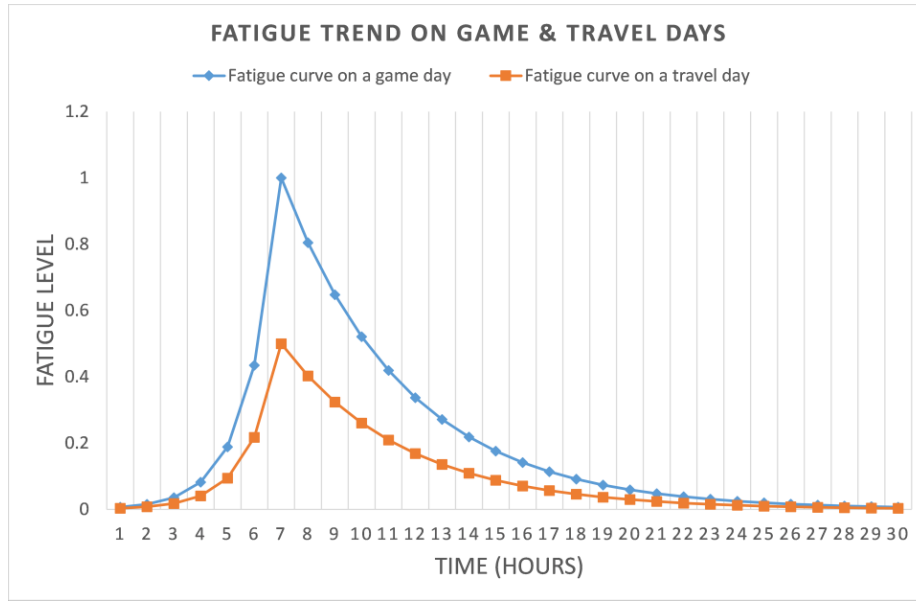

**Figure 2:** General fatigue trends on play and travel days

*4.3.1 Maximum Fatigue*

MS was the franchise that suffered the maximum fatigue of 1.3456 during the fifth edition of the PSL, followed by QG, IU, PZ, LQ, and then KK at the last with 1.1541. On the other hand, the average peak fatigue by each franchise in alphabetical order is 1.13, 1.10, 1.08, 1.12, 1.11, 1.14; which favors LQ the most, followed by KK, PZ, MS, IU, and QG in order.

*4.3.2 Prior Average Fatigue*

Average fatigue value at the start of each game and travel for IU, KK, LQ, MS, PZ and QG are recorded as 0.166, 0.118, 0.105, 0.195, 0.135, and 0.175 respectively. It is once again evident that this schedule favored LQ the most, which is around 85% better than the least favored MS.

*4.3.3 Separation Between Successive Games*

The average separation between successive game peaks is recorded as 53.3 hours for PZ, 61.3 hours each for IU, KK and MS, 62 hours for LQ, and finally 63.3 hours for QG. These numbers show that PZ had the most constrained schedule, whereas QG had the most relaxed, with plenty of days to rest before each game. However, the average does not always reflect the reality, since, in this case for example, poor distribution of rest days – very few before some games and unnecessarily many before others – may yield similar averages, but vastly varying fatigue levels. In order to present a relatively simplified view of the matter, we have computed utilization of the available rest days per game that each team plays; this is presented in Fig. 6. Note the linear plot labelled *Trend*, which presents the ideal rest days utilization per game. The team having a curve closest to the ideal situation must be considered the most favored one. For this purpose, we utilize RMSE as the competitive metric between the franchises' schedules. Tab. 5 summarizes the computed RMSE values for each team. Clearly, KK's schedule most closely fits the ideal rest days utilization, which is around 55% better than the worst schedule of IU.

*4.4 Conclusive Remarks on Fairness Using a Weighting Criterion*

Based on the above statistics, the readers may formulate a mechanism on their own to nominate the most favored side by the PSL-5 schedule. For reference, we have devised a simple weighting criterion that helps us in pointing out the most favored team by the schedule as follows:

i. The sides having the smallest and the largest RMSE get 5 and 0 points respectively; the rest are

adjusted in between.
ii. The sides having the lowest and the highest average fatigue values at the start of each game get 5 and 0 points respectively; the rest are adjusted in between.
iii. The sides having the lowest and the highest average of maximum fatigue values get 5 and 0 points respectively; the rest are adjusted in between.
iv. The points ($P$) for traveling ($T$) by each side ($i$), are computed using $P_{Ti} = |T_i - \max(T_j)| \forall j \in \{1 \ldots 6\}$. For example, MS traveled six times, so $T_{MS} = 6$, which is also the maximum amongst all

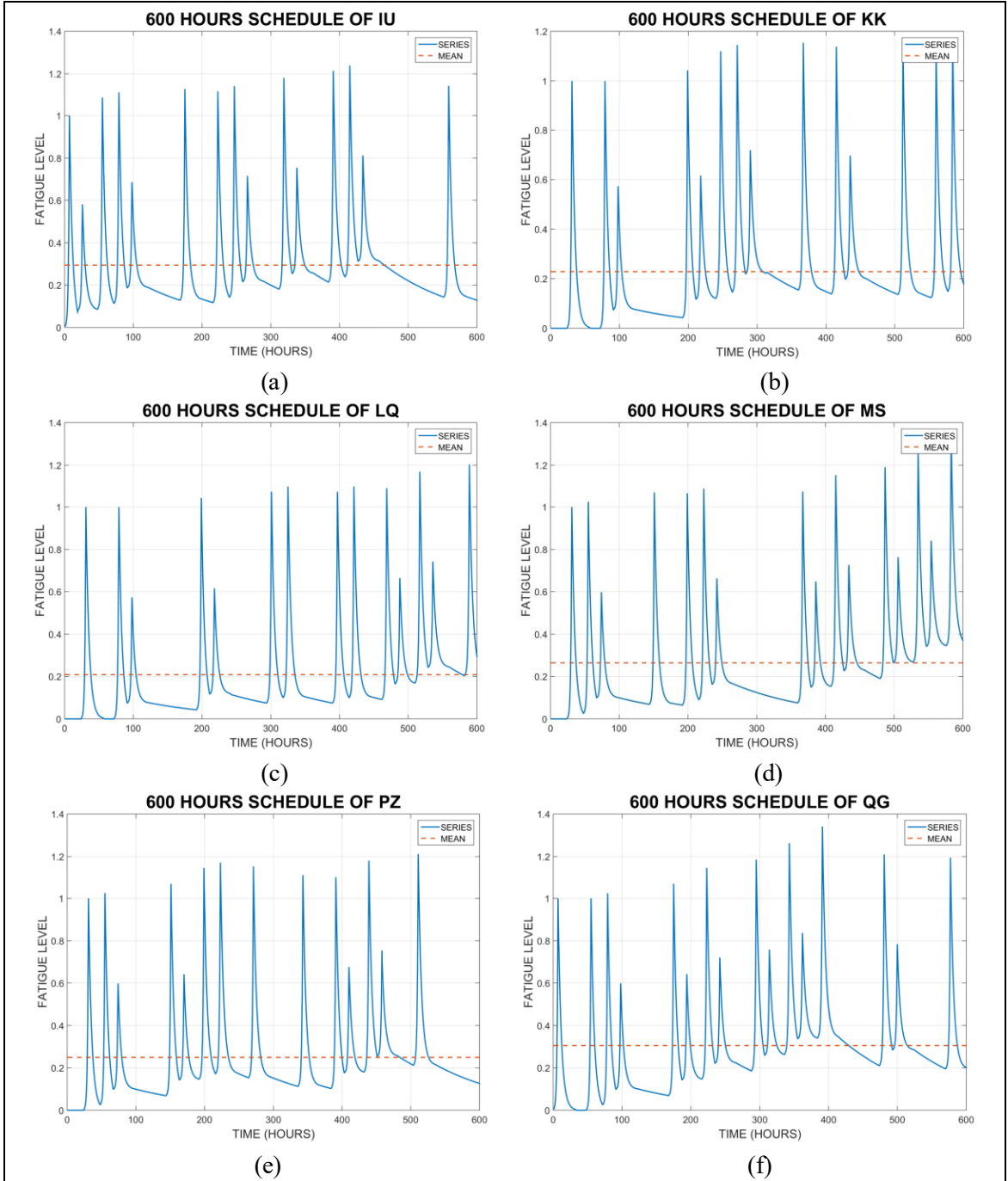

**Figure 3:** Fatigue trend for each franchise during the 600 hours tournament span

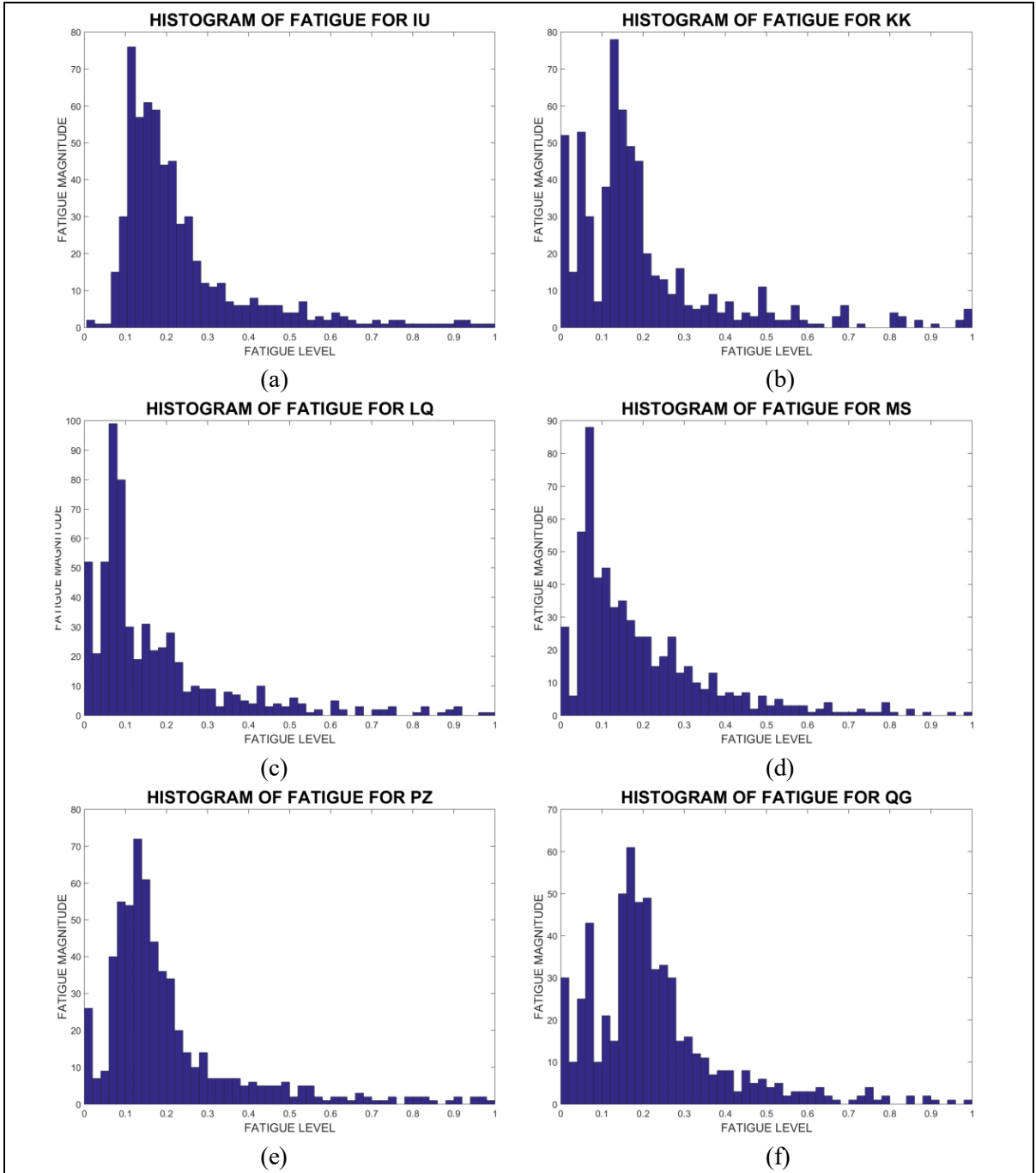

**Figure 4:** Fatigue histogram for each franchise during the 600 hours tournament span

the sides, leading to $P_{TMS} = 6$.

    v.    Each rest day earns one point.

    vi.    Each home game earns one point.

Note that items 1 to 5 above add to overall fatigue, while the sixth one offers a different kind of advantage. Nevertheless, we add all these up to determine the most favored side - the result in presented in Tab. 6. As far as the fatigue is concerned, KK, having maximum points (sum of columns 1 to 5), seems to have enjoyed this schedule the most, followed by LQ, and MS at the bottom. Considering the home advantage (column 6), however, LQ appears to be the most favored franchise during the previous edition of the PSL.

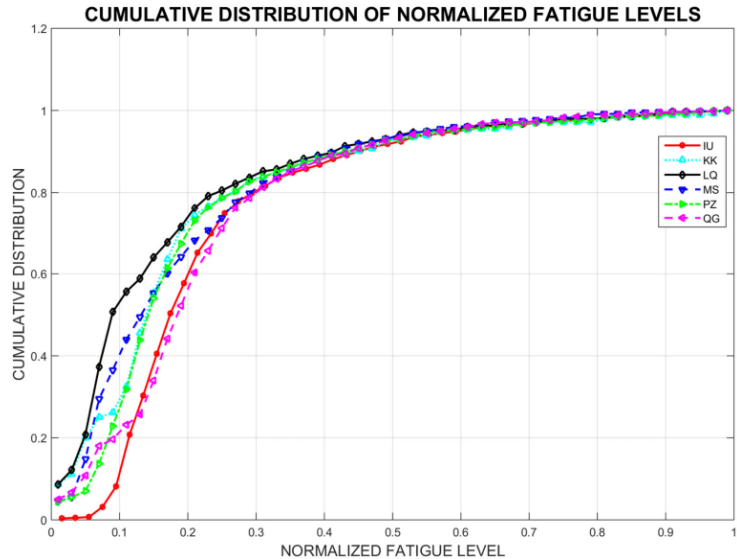

**Figure 5:** Cumulative Distribution of Fatigue Levels

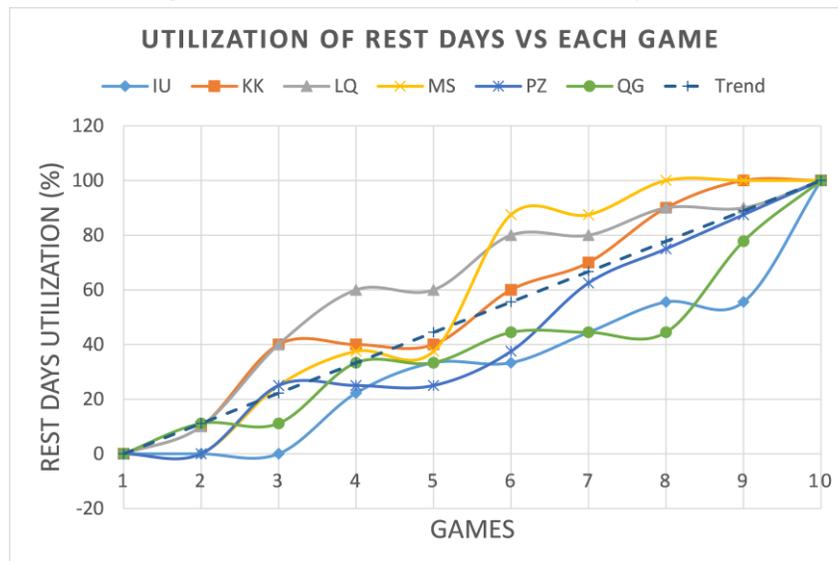

Figure 6: **Distribution of rest days across the ten round games**

*4.5 Outstanding Research Objectives*

There is unfairness within the schedule at two levels discussed below, addressing which defines two critical outstanding research objectives.

i. Fairness with respect to competitive parameters: Fig. 1 presented the first level at which fairness is required. There must be uniformity in each competitive parameter among the franchises.

ii. Fairness with respect to distribution of rest days: Fig. 6 has presented RMSE as the other level of unfairness. Each team deserves to have a schedule with minimal RMSE. In this case, we must ensure uniformity with respect to RMSE for each side as well.

**Table 5:** RMSE value for each franchise

| Team | IU | KK | LQ | MS | PZ | QG |
|---|---|---|---|---|---|---|
| RMSE | 18.59 | 8.27 | 14.81 | 15.06 | 9.65 | 14.48 |

Table 6: A concluding quantitative comparison of unfairness

| Items | 1 | 2 | 3 | 4 | 5 | 6 | Fatigue Points | Total Points |
|---|---|---|---|---|---|---|---|---|
| IU | 0 | 2 | 1 | 1 | 9 | 5 | 13 | 18 |
| KK | 5 | 4 | 4 | 2 | 10 | 5 | 25 | 30 |
| LQ | 2 | 5 | 5 | 2 | 10 | 8 | 24 | 32 |
| MS | 1 | 0 | 2 | 0 | 8 | 3 | 11 | 14 |
| PZ | 4 | 3 | 3 | 2 | 8 | 3 | 20 | 23 |
| QG | 3 | 1 | 0 | 0 | 9 | 2 | 13 | 15 |

In the sequel of this article (part 2), we will develop a set of novel constraints specifically formulated to address the aforementioned unfairness with the PSL-5 schedule, present the case as a multiobjective optimization problem [26], and propose a framework based on integer programming to yield optimal schedules. The same may be adopted for other leagues/ sports as well, however, considering their strictly time bound nature, our proposed framework will particularly suit the T20 cricket leagues played worldwide.

## 5 Conclusion

We have considered one of the schedules used for a recent edition of the Pakistan Super League and demonstrated the unfairness in various terms within it. We have exploited various competitive parameters to carry out a comparison between the schedules of the six participating teams. The parameters include the total number of days to complete their round games, the number of rest days, the number of home and away games, and the number of games that each side plays successively. The last one is the most crucial, since it defines the distribution of rest days across the entire schedule. We have used root mean squared error between an ideal and the actual schedules, and cumulative distribution function to point out the sides that this schedule favored the most. We have advocated for optimal tournament scheduling, which must ensure fair distribution of rest days to avoid sportsmen burnout.

A natural follow-up work, which is lacking in this manuscript, is to propose a novel scheduling framework that should rectify the unfairness concerns pointed out in this work. However, in part 2 of the same manuscript, we will present novel constraints to address the issues highlighted in this part of the article, and we propose a framework for optimizing sports tournament schedules. While, our framework will be applicable to any sports event, we specifically focus on the T20 cricketing events, considering its strictly time bound nature.

**Funding Statement:** The authors are grateful to the Deanship of Scientific Research at King Saud University, Saudi Arabia for funding this work through the Vice Deanship of Scientific Research Chairs: Chair of Pervasive and Mobile Computing.

**Conflicts of Interest:** The authors declare that they have no conflicts of interest to report regarding the present study.